\documentclass[letterpaperpaper,11pt]{article}
\usepackage{latexsym,amsmath,amssymb}
\usepackage{epsfig,graphics}
\advance\topmargin -1cm
\advance\oddsidemargin -1.5cm
\advance\evensidemargin -1.5cm
\newcommand{\cqd}{\hfill\rule{2mm}{2mm}}
\def\car{\mbox{\rm{1\hspace{-0.10 cm }I}}}

\newtheorem{prop}{Proposition}[section]
\newtheorem{teor}{Theorem}[section]

\newtheorem{remark}{Remark}[section]
\newtheorem{ejemplo}{Example}[section]
\newtheorem{lema}{Lemma}[section]
\begin{document}
\title{ Asymptotic properties of the maximum likelihood
estimator for nonlinear AR processes with markov-switching}

\author{Luis-Angel Rodr\'{\i}guez\\
{\small Dpto. de Matem\'aticas, FACYT, Universidad de Carabobo, Venezuela}\\
{\small CIMFAV, Facultad de Ingenier\'{\i}a, Universidad de Valpara\'{\i}so, Chile.}\\
{\small Corresponding: larodri@uc.edu.ve}}

\maketitle

\begin{abstract}
In this note, we propose a new approach for the proof of the consistency and normality of the maximum likelihood estimator
for nonlinear AR processes with markov-switching under the assumptions of uniform exponential forgetting of the prediction filter and
$\alpha$-mixing pro\-per\-ty. We show that in the linear and Gaussian case our assumptions are fully satisfied.\\
\noindent{{\bf Keywords:} Nonlinear autoregressive process, Markov switching
asymptotic normality,consistency, hidden Markov chain.}\\
\noindent{{\bf MSC:} Primary: 60G17, Secondary:62G07}
\end{abstract}

Switching autoregressive processes with Markov re\-gi\-me can be considered as a combination of hidden Markov models (HMM) and threshold
regression models. They have been introduced in an econometric
context by Goldfeld and Quandt (1973) \cite{Goldfeld} and they have
become quite popular in the literature  since Hamilton (1989)
\cite{Hamilton} employed them in the analysis of the the rate of growth of USA GNP series
for two regimes: one of contraction and another of expansion. This family of models  describes 
the evolution of a time series subject to discrete shifts and
the transition is controlled by a HMM. \\

We consider a nonlinear AR process with markov-switching (abbreviated MS-NAR) $\{Y_n\}_{n\geq0}$ defined 
for integers $n\geq 1$ by
\begin{equation}
    \label{modelo}
    Y_n=r(Y_{n-1},\theta_{X_n})+e_n,\ Y_n\in\mathbb{R}.
\end{equation}
Here the process $\{e_n\}_{n\geq1}$ are i.i.d. random variables and 
the sequence $\{X_n\}_{n\geq1}$ is an homogeneous Markov chain with state space $\{1,\ldots,m\}$. \\

Let $\mathcal{F}=\{r(\cdot,\theta):\theta\in\Theta\}$  
a family of real valued functions defined on $\mathbb{R}^{m+1}$,
indexed by a parameter $\theta=(\theta_1,\ldots,\theta_m)\in\Theta$ 
and $\Theta$ is a compact set of $\mathbb{R}^{m}$. 
We denote by $A$ the probability transition matrix of 
the Markov chain  $\{X_n\}_{n\geq1}$, i.e. $A=[a_{ij}]$, with $a_{ij}=\mathbb{P}(X_n=j|X_{n-1}=i)$. 
The parameter space is the set
$$
\Psi= \left\{\psi= (\theta,A):\theta\in\Theta, a_{ij}\in[0,1]\ \mbox{and} \sum_{j=1}^m a_{ij}=1\right\}.
$$

We assume that the variable $Y_0$, the Markov
chain $\{X_n\}_{n\geq1}$ and the sequence $\{e_n\}_{n\geq1}$ are mutually
independent. The process $\{X_n\}$, called regime, is not observable and inference 
has to be carried out in terms of the observable process $\{Y_n\}$. \\

The consistency of the maximum likelihood estimator for the parameter $\psi$ in the
MS-NAR model is given in Krishnarmurthy and Ryden (1998) \cite{kris-Ryden}, 
while the consistency and asymptotic normality are proved in a more general context in
the work of Douc {\it et al.} (2004) \cite{douc}. In the section 2 we prove the consistency and asymptotic normality of the maximum likelihood
estimator for functional AR processes with markov-switching under the assumptions of exponential uniform forgetting property for prediction filter and
an $\alpha$-mixing property. 

\section{General properties for MS-NAR model}
In this section we review the key properties of the
MS-NAR model that we need for proving our results.

\subsection{Stability and existence of moments}

The study of the stability of the model MS-NAR is relatively complex. In this section we recall known
results about the stability of this model given by Yao and Attali \cite{Yao}. Our aim is to resume the sufficient conditions 
which ensure the existence and the uniqueness of a stationary ergodic solution for the model, as well as the existence of moments 
of order $s\geq 1$ of the respective stationary distribution.  
 
\begin{enumerate}
\item[\bf E1] The Markov chain $\{X_n\}_{n\geq 1}$ is positive recurrent. Hence, it has an invariant distribution that we 
denote by $\mu=(\mu_1,\ldots,\mu_m)$.

\item[\bf E2] The functions $y\to r(y,\theta_i)$, for $i=1,...,m$, are continuous. 

\item[\bf E3] There exist positive constants $\rho_i, b_i$, $i=1,...,m$,  such that for $y\in\mathbb{R}$, the following
inequality holds
$$
|r(y,\theta_i)|\leq \rho_i|y|+b_i.
$$

\item[\bf E4] $\gamma=\sum_{i=1}^m\log \rho_{i}\mu_i<0$.

\item[\bf E5] $\mathbb{E}(|e_1|^s)<\infty$, for some $s\geq 1$.

\item[\bf E6] The sequence $\{e_n\}_{n\geq1}$  of random variables admits a common density probability function $\Phi$ 
with respect to the Lebesgue measure.

\item[\bf E7] There exist $b>0$ and $\mathrm{C}$ a compact set of $\mathbb{R}$ such that $\inf_{e\in\mathrm{C}}\Phi(e)>b$.
\end{enumerate}

Condition E1 implies that $\{(Y_n,X_n)\}_{n\geq1}$
with states space  $\mathbb{R}\times\{1,\ldots,m\}$ is a  Markov process. Under condition E2 this Markov process 
is a Feller chain and it is a strong Feller chain if in addition the condition E6 holds.\\

The model is called sublinear if conditions E2 and E3 hold. For the sublinear MS-NAR model, 
Yao and Attali \cite{Yao} proved the following result.

\begin{prop} Consider a sublinear MS-NAR $\{Y_n\}_{n\geq0}$. Under assumptions 
E1-E7, we have that 
\begin{itemize}
\item[i)] There exists a unique stationary geometric ergodic solution.

\item[ii)] If the spectral radius of the matrix $Q_s=\left( \rho_j^s \, a_{ij}\right)_{i,j=1\ldots m}$ is strictly less than 1, 
with $s$ given in  E5, then $\mathbb{E}(|Y_n|^s)<\infty$.
\end{itemize}
\end{prop}

\begin{remark} The Markov chain is stable under the moment condition $s\geq 1$, but for
the asymptotic properties of the MLE it will be necessary to assume $ s>2$.
\end{remark}

Now we introduce some notations:

\begin{itemize}
\item $V_{1:n}$ stands for the random vector $(V_1,\ldots,V_n)$, and by $v_{1:n}=(v_1,\ldots,v_n)$ we mean
a realization of the respective random vector.
\item The symbol $\car_B(x)$ denotes the indicator function of set $B$, which assigns the value $1$ if $x\in B$ and  $0$ otherwise.
\item $p(V_{1:n}=v_{1:n})$ denotes the density distribution of random vector $V_{1:n}$ evaluated at $v_{1:n}$.
\end{itemize}


We consider the following assumption :

\begin{enumerate}
\item[\bf D1] The random variable $Y_0$ admits a density function $p(Y_0=y_0)$ with respect to Lebesgue measure.
\end{enumerate}

Under conditions D1 and E6, the random vector $(Y_{0:n},X_{1:n})$ admits the probability density
$p(Y_{0:n}=y_{0:n},X_{1:n}=x_{1:n})$ equal to
$$
\Phi(y_n-r(y_{n-1},\theta_{x_n}))\cdots\Phi(y_{1}-r(y_0,\theta_{x_n}))a_{x_{n-1}x_n}\cdots
a_{x_{1}x_2}\mu_{x_1}p(Y_0=y_0),
$$
with respect to the product measure $\lambda\otimes\mu_c$, where $\lambda$ and $\mu_c$
denote Lebesgue and counting measures respectively. For a proof of this result see
Ferm\'{\i}n {\it et al} \cite{luis3}.

\subsection{Strong mixing}
A strictly stationary stochastic process $Y =\{Y_n\}_{n\in\mathbb{Z}}$ is called strongly mixing, if 
\begin{equation}
\alpha_n :=\sup\{|\mathbb{P}(A\cap B)- \mathbb{P}(A)\mathbb{P}(B)|: A\in \mathcal{M}_{-\infty}^0, B\in \mathcal{M}_n^{\infty}\} \to 0, \quad \mbox{ as } n\to \infty ,
\end{equation}
where $\mathcal{M}_a^{b}$, with $a,b \in \overline{\mathbb{Z}} $, is the $\sigma$-algebra generated by $\{Y_k\}_{k= a:b}$, and is 
absolutely regular mixing if 
\begin{equation}
\beta_n :=\mathbb{E}\left( ess \sup\{\mathbb{P}(B| \mathcal{M}_{-\infty}^0)-\mathbb{P}(B) :  B\in \mathcal{M}_n^{\infty}\} \right) \to 0, 
\quad \mbox{ as } n\to \infty.
\end{equation}

The values $\alpha_n$ and $\beta_n$ are called $\alpha$-mixing and $\beta$-mixing coefficients respectively. 
For properties and examples of processes under mixing assumptions, see Doukhan \cite{dmixing}. 
In general, we have the inequality $2\alpha_n \leq \beta_n\leq 1$.\\

Note that the $\alpha$-mixing coefficients can be rewritten as:
\begin{equation}
\alpha_n :=\sup\{|cov(\phi,\xi)|: 0\leq \phi, \xi \leq 1, \phi \in \mathcal{M}_{-\infty}^0, \xi\in \mathcal{M}_n^{\infty}\}. 
\end{equation} 
In the case of a strictly stationary Markov process $X$, with state space $(E,\mathcal{B})$, kernel probability transition $A$ 
and invariant probability measure $\mu$, the $\beta$-mixing coefficients take the following form (see Doukhan \cite{dmixing}, section 2.4):
\begin{equation}
\beta_n :=\mathbb{E}\left( \sup\{|A^{(n)}(X,B)-\mu(B)|: B\in \mathcal{B}\} \right).
\end{equation}

\begin{lema}
\label{ModeloesMixing}
Under conditions E1-E7 the process MS-NAR  is $\alpha$-mixing with $\alpha$-mixing coefficients 
decreasing geometrically.
\end{lema}

\noindent{\bf Proof:} For the proof of this lemma see Ferm\'{\i}n {\it et al} \cite{luis3}. 
\cqd\\

\begin{ejemplo} (Linear autoregressive with Markov switching (MS-AR)  nonmixing)
\end{ejemplo}
In the case where $r(y,(b_i,\rho_i)^t)=\rho_{i}y+b_{i}$, 
the model is a MS-AR and it is defined by:
\begin{equation}
\label{modelolineal}
         Y_{n}=\rho_{X_n}Y_{n-1}+b_{X_n}+e_n.
\end{equation}
For each $1\leq i\leq m$, we denote $\theta_i=(b_i,\rho_i)^{t}$ and
 $$
\theta=\left(
\begin{array}{cccc}
b_1 &b_2&\cdots &b_m\\
\rho_1 &\rho_2&\cdots & \rho_m
\end{array}\right).
$$
More specifically consider the process MS-AR  with $\theta_i=(0,\rho_i)^t$
for all $i=1,\ldots,m$ and such that the random variable $e_1$ follows 
a Bernoulli distribution with parameter $q$ and $Y_0=0$. In this case,
we have
$$
Y_n=\sum_{k=0}^{n-1}\rho_{X_k}\cdots\rho_{X_1}e_{k+1},
$$ 
and we adopt the convention that $\rho_{X_k}\cdots\rho_{X_1}=1$ for $k=0$. This process is non $\alpha$-mixing.  In fact, according to D. Andrews \cite{Andrews} if 
$0<\rho_i\leq 1/2$, for  $t\in \mathbb{N}$ there exist some sets  $A\in\mathcal{M}_{-\infty}^0$, 
$B_t\in \mathcal{M}_n^{\infty}$,  with $\mathbb{P}(A)>0$,
$\mathbb{P}(B_s)\leq c$  for some constant $c<1$ such that $\mathbb{P}(B_t|A)=1$,  therefore
$$
\alpha_t(Y)\geq \mathbb{P}(A\cap B_t)-\mathbb{P}(A)\mathbb{P}(B_t)=\mathbb{P}(A)(\mathbb{P}(B_t|A)-\mathbb{P}(B_t))\geq \mathbb{P}(A)(1-c).
$$
This implies that $\alpha_t(Y)$ does not tend to $0$ as $t\to\infty$ 
and so  $Y$ is a non $\alpha$-mixing process.

\begin{lema}
\label{LeyesMixing}
Under conditions E1-E7, the MS-NAR process $\{Y_n\}_{n\geq0}$ satisfies,
\begin{itemize}
\item[i)] For all function $\varphi$ such that $\mathbb{E}(\varphi(Y_{k}))<\infty$, 
we have the strong law of large numbers,
$$
\frac{1}{n}\sum_{k=1}^n\varphi(Y_{k})\to \mathbb{E}(\varphi(Y_1)),\ a.s.
$$
\item[ii)] Suppose that $\mathbb{E}(\varphi(Y_1))=0$, $\mathbb{E}|\varphi(Y_1)|^{s}<\infty$, for some $s>2$.
Then
$
\Gamma=\mathbb{E}(\varphi(Y_1)^2)+2\sum_{k=1}^{\infty}k+\mathbb{E}(\varphi(Y_1)\varphi(Y_k))<\infty 
$
and if $\Gamma\not=0$, 
$$
\frac{1}{\sqrt{n}}\sum_{k=1}^n\varphi(Y_{k})\to \mathcal{N}(0,\Gamma), 
$$
for $n\to\infty$, in distribution.
\end{itemize}
\end{lema}

\noindent{\bf Proof:} i)  This result is a direct consequence of the Collolary 3.1 in Rio \cite{Rio1}. \\

\noindent ii) Let $U_k=\varphi(Y_{k})$, then $\{U_k\}_{k\geq0}$ is a strictly stationary sequence and is strongly $\alpha$-mixing, with
$\mathbb{E}(|U_k|^r)<\infty$, for $r>2$. For $\alpha^{-1}(u)=\inf\{k\in\mathbb{N}:\ \alpha_k\leq u \}$, we 
have to prove
\begin{equation}
  \label{DM}
  \int_{0}^1\alpha^{-1}(u)Q^2(u)du<\infty 
\end{equation}
where  $Q$ is the associate quantile function of the process $\{U_k\}$.
The condition \eqref{DM} is implied by
\begin{equation}
\label{condicionTCL}
\sum_{i\geq0}(i+1)^{\frac{2}{r-2}}\alpha_i<\infty,
\end{equation}
and in our case this is valid, since from geometric $\alpha$-mixing property 
exist $0<\zeta<1$ such that $\alpha_i\leq C \zeta^i$, we have
$$
\sum_{i\geq0}(i+1)^{\frac{2}{r-2}}\alpha_i\leq C\sum_{i\geq0}(i+1)^{\frac{2}{r-2}}\zeta^i<\infty.
$$
Thus, we can apply Theorem 4.2. in E. R\'{\i}o \cite{Rio1}, obtaining 
that $\sqrt{n}U_n$ converges in distribution to $\mathcal{N}(0,\Gamma)$.
\cqd
\section{Maximum likelihood estimation}

Using $p_\psi$ as a generic simbol for densities and distributions parameterized for $\psi$. We defined
the conditional $\log$-likelihood as $l_n(\psi)=\log p_\psi(Y_{1:n}|Y_0)$ and
we can expressed as
$$
l_n(\psi)=\sum_{k=1}^n \log p_\psi(Y_k|Y_{0:k-1}).
$$

We denote by $\psi^*$ the true parameter wich is consider as fixed.
A maximum likelihood estimator (MLE) is defined by
$$
\hat{\psi}_n=\arg\max_{\psi} l_n(\psi).
$$
The MLE is consistent if $\hat{\psi}_n\to\psi^*$ as $n\to\infty$ a.s.\\ 

The techniques standard used to prove consistency follows the steps:
\begin{enumerate}
\item To show that there exists a continuous deterministic function $l(\psi)$ such that 
 $$
\lim_{n\to\infty} \frac{1}{n}l_n(\psi)=l(\psi)\ \mbox{a.s.}$$  
\item To show that $l(\psi)$ a.s. has a unique maximum at $\psi=\psi^*$.
\item To conclude that $\hat{\psi}_n=\arg\max_{\psi} n^{-1}l_n(\psi)\to\arg\max_{\psi} l(\psi)=\psi^*$.
\end{enumerate}
\vspace{0.5cm}

For MS-NAR processes a strong law of large numbers of the $log$-likelihood is 
obtained in Rynkiewicz \cite{ryn}, Krish\-na\-mur\-thy \cite{Kris} 
using an additive function of the extended Markov chain $(Y_n,X_n,\mathbb{P}_\psi(X_k|Y_{0:n}))$.
In Douc {\it et. al.} \cite{douc} the law of large numbers of the $log$-likelihood
follow from uniform exponential forgetting of the initial
distribution for prediction filter.\\
 
In this work, following the approach of consistency proof of
Handel, chapter 7 in \cite{Ramon}, for HMM,
and joined to the $\alpha$-mixing property we obtain a new proof
of the consistency for the MLE.\\

The following lemma shows that we can express $p_\psi(Y_{k}|Y_{0:k-1})$
as a functional of the prediction filter $\mathbb{P}_\psi(X_k|Y_{0:n})$. 

\begin{lema}
\label{filtro}
Let $\delta=\inf_{i,j=1:m} a_{ij}$. Define,
$$D_{k,l}^\psi=\log \int\int p_\psi(Y_k|Y_{0:k-1},x_k) a_{x_{k-1},x_{k}}\mathbb{P}(x_{k-1}|Y_{l:{k-1}})\mu_c(dx_k)\mu_c(dx_{k-1}),$$
for $0<l<k$. Under assumptions E1 and E7 then $|D_{k,l}^\psi-D_{k,0}^\psi|\leq2\delta^{-1}(1-\delta)^{k-1-l}$.
\end{lema}

\noindent {\bf Proof:} First, we bound from below the quantities $\exp(D_{k,0}^\psi)$ y $\exp(D_{k,l}^\psi)$,
by the Fubini Theorem we have
$$
\exp(D_{k,0}^\psi)\geq\delta\int p_\psi(Y_k|y_{0:k-1},x_k) \mu_c(dx_k)
$$
and the same for $\exp(D_{k,l}^\psi)$, thus 
$$
\min(\exp(D_{k,0}^\psi),\exp(D_{k,l}^\psi))\geq \delta\int p_\psi(Y_k|Y_{0:k-1},x_k) \mu_c(dx_k).
$$

Using inequality $|\log x- \log y|\leq|x-y|/\min(x,y)$, we estimate
\begin{eqnarray*}
\lefteqn{|D_{k,l}^\psi-D_{k,0}^\psi|}\\
&\leq& \frac{\int\int p_\psi(Y_k|Y_{0:k-1},x_k) a_{x_{k-1},x_{k}}(\mathbb{P}_\psi(x_{k-1}|Y_{l:{k-1}})-
\mathbb{P}_\psi(x_{k-1}|Y_{0:{k-1}}))\mu_c(dx_k)\mu_c(dx_{k-1})}{\delta\int p_\psi(Y_k|Y_{0:k-1},x_k) \mu_c(dx_k)}\\
&\leq& \frac{1}{\delta}
\|\mathbb{P}_\psi(X_{k-1}\in\cdot|Y_{l:{k-1}})-\mathbb{P}_\psi(X_{k-1}\in\cdot|y_{0:{k-1}})\|_{VT}
\end{eqnarray*}
Applying the Proposition 4.3.26 (iii) in Cappe {\it et. al. } \cite{LibCappe}, p\'ag 109,
$$
\|\mathbb{P}_\psi(X_k\in\cdot|Y_{l:k})-\mathbb{P}_\psi(X_k\in\cdot|Y_{0:k})\|_{VT}\leq2(1-\delta)^{k-l}
$$
We conclude that $|D_{k,l}^\psi-D_{k,0}^\psi|\leq2\delta^{-1}(1-\delta)^{k-1-l}$.\cqd\\

This lemma shows that the quantity $D_{k,0}^\psi$, which
depends on the observations $Y_{0:k}$, can be approximated by
$D_{k,l}^\psi$ which is a function of only a fixed number of observations  $Y_{l:k}$.

\begin{prop}
Under assumptions E1-E7, suppose $\Psi$ is a compact set and
the condition $\mathbb{E}_{\psi^*}(p_{\psi}(Y_k|Y_{k-1},i))<\infty$, for $i=1,\ldots,m$. Then
$l_n(\psi)$ is a continuous function and $l(\psi)=\lim_{n\to\infty}{n}^{-1}l_n(\psi)$ exist a.s 
for each $\psi\in\Psi$.
\end{prop}

\noindent{\bf Proof:} The proof is done in two steps. First, we have
$$
l(\psi)=\lim_{k\to\infty}\mathbb{E}_{\psi^*}(D_k^{\psi}), \mbox{exists for every $\psi\in\Psi$.} 
$$

Second, we show
$$
\frac{1}{n}\sum_{k=1}^n(D_k^\psi-\mathbb{E}_{\psi^*}(D_k^{\psi}))\to 0\ a.s.
$$
then conclude 
$$
\frac{1}{n}l_n(\psi)=\frac{1}{n}\sum_{k=1}^n(D_k^\psi-\mathbb{E}_{\psi^*}(D_k^{\psi}))+
\frac{1}{n}\sum_{k=1}^n\mathbb{E}_{\psi^*}(D_k^\psi)\to l(\psi)\ a.s.
$$
\noindent{\bf Step 1}.  Let $\Delta_k=\mathbb{E}_{\psi^*}(D_{k,l}^{\psi})$ by Lemma \ref{filtro},
$$
|\Delta_{m+n}-\Delta_m|=|\mathbb{E}_{\psi^*}(D_{m+n,0}^{\psi}-\mathbb{E}_{\psi^*}(D_{m+n,m}^{\psi}))|\leq2\delta^{-1}(1-\delta)^{m-1}
$$
hence $\sup_{n}|\Delta_{m+n}-\Delta_m|\to0$ as $m\to\infty$, i.e.,
$\{\Delta_k\}$ is a Cauchy sequence and  therefore convergent. By Ces\`{a}ro’s theorem
$\mathbb{E}_{\psi^*}(l_n(\psi))=n^{-1}(\sum_{k=0}^{n-1}\Delta_k)$ 
also converges.\\

\noindent{\bf Step 2}. According to Proposition 1.1 the sequence $\{D_k^\psi\}_{k\geq1}$ is $\alpha$-mixing
with geometric coefficients $\alpha_k$. We demostrate that  
$\mathbb{E}(|D_k^\psi|)<\infty$, in fact
$$
\int\int p_\psi(Y_k|Y_{0:k-1},x_k) a_{x_{k-1},x_{k}}\mathbb{P}_\psi(x_{k-1}|Y_{l:{k-1}})\mu_c(dx_k)\mu_c(dx_{k-1})
\leq \int p_\psi(Y_k|Y_{0:k-1},x_k) \mu_c(dx_k)
$$
and 
$$
\int p_\psi(Y_k|Y_{0:k-1},x_k) \mu_c(dx_k)=\sum_{i=1}^m p_\psi(Y_k|Y_{k-1},i)\mu_i\leq m\max_{i=1:m}\{p_\psi(Y_k|Y_{k-1},i)\mu_i\} 
$$
under assumption $\mathbb{E}(p_{\psi^*}(Y_k|Y_{k-1},i))<\infty$, then
$
\mathbb{E}(|D_k^\psi|)<\infty
$
and by Lemma \ref{LeyesMixing}, i) we obtain
$$
\frac{1}{n}\sum_{k=1}^n(D_k^\psi-\mathbb{E}_{\psi^*}(D_k^{\psi}))\to 0\ a.s.
$$
\cqd\\

We prove the validity of step three under  uniform convergence,
$\sup_{\psi\in\Psi}|l_n(\psi)-l(\psi)|\to 0$.

\begin{lema}
Suppose $\Psi$ is a compact set. Let $l_n:\Psi\to\mathbb{R}$ be a sequence of con\-ti\-nuous
functions that converges uniformly to a function $l:\Psi\to\mathbb{R}$. Then
$$
\hat{\psi}_n=\arg\max_{\psi} l_n(\psi)\to\arg\max_{\psi} l(\psi)
$$
\end{lema}
\noindent {\bf Proof:} As a continuous function on a compact space attains its maximum, we
can find a $\psi_n\in\arg\max_{\psi} l_n(\psi)$ for all $n$.
Which show using an argument that goes to Wald (1949) that 
\begin{equation}
\label{star1}
\lim_{n\to\infty} l(\psi_n)=\sup_{\psi\in\Psi} l(\psi). 
\end{equation}

Suppose that the sequence $\{\psi_n\}$ does not converge to the set 
$\{\tilde{\psi}:\ l(\tilde{\psi})=\max_{\psi\in\Psi}l(\psi)\}$. By compactness there
exists a subsequence $\{\psi_n'\}\subset\{\psi_n\}$ which converges 
to $\psi'\not\in\{\tilde{\psi}:\ l(\tilde{\psi})=\max_{\psi\in\Psi}l(\psi)\}$. 
But $l(\psi)$ is continuous, so  $l(\psi_n')\to l(\psi')<\sup_{\psi\in\Psi} l(\psi)$ 
and according to \eqref{star1}, this is a contradiction.
\cqd

\begin{teor}
Suppose $\Psi$ is a compact set. Assume that
\begin{enumerate}
 \item $\psi=\psi^*$ iff $\mathbb{P}_\psi=\mathbb{P}_{\psi^*}$.
 \item For all $i,j\in\{1,\ldots,m\}$ and all $y,y'\in\mathbb{R}\times\mathbb{R}$ 
 the functions $\psi\to a_{ij}$ and $\psi\to p_\psi(Y_1=y|Y_0=y',X_1=i)$ are continous. 
 \item There is a $c<\infty$ such that 
 $|D_k^\psi-D_k^{\psi'}|\leq c\|\psi-\psi'\|$
 for all $k>1$
\end{enumerate}
Then the maximum likelihood estimate $\hat{\psi_n}$ is consistent.
\end{teor}

\noindent{\bf Proof:} By Theorem 7.5 in Handel \cite{Ramon}, the Lipschitz condition 3. and 
compactness implies that the sequence $l_n\to l$ a.s uniformly. According to Lemma 2.1
$$
\hat{\psi}_n\to \psi_{*}=\arg\max_{\psi} l(\psi),
$$
and this value is unique under identifiability.\cqd\\ 

In the Gaussian and linear case we can prove directly identifiability and equicontinuity.
This allows us obtain the consistency of the MLE  
without assuming a condition of Lipschitz for the parameters.

\begin{ejemplo}(MS-AR gaussian linear)
\end{ejemplo}

Let the model defined by \eqref{modelolineal}. Let
$\{e_n\}$ are gaussian i.i.d. random variables.  
Our goal in this example is check that the conditions for
consistency apply in this case. In fact, if
we assume that for the true model
$\Psi^{*}$ the vector components
$\{(\alpha_i,b_i,\sigma_i)\}_{i=1}^{m}$ are different; thus, for every
$n$, there exists a point $Y_{n-1}\in\mathbb{R}$ such that
$\{(\alpha_iY_{n-1}+b_i,\sigma_i)\}_{i=1}^{m}$ are different.
Therefore, in agreement with Remark 2.10 of Krishnamurthy and Yin \cite{Kris}
the model is identifiable in the following sense:
If  $K$ stands for the Kullback-Leibler divergence $K(\psi,\psi_{*})=0$ then, 
$\psi=\psi_{*}$, which proves the identifiability. 
On the another hand, the Lemma 4.1 in \cite{luis1} follows 
that $\frac{1}{n}\log p_\psi(Y_1^n|Y_0=y_0)$ is an equicontinuos sequence a.s-$\mathbb{P}_{\psi_{*}}$.
We conclude that in this case the MLE is consistent.\\

There is a standard technique for prove asymptotic normality of
maximum likelihood estimates. The idea is that the first derivatives of a
smooth function must vanish at its maximum. If we
expand in Taylor series the likelihood gradient around  $\psi^*$,
we can write
$$
0=\nabla_\psi l_n(\hat{\psi}_n)=\nabla_\psi l_n(\psi^*)+\nabla_\psi^2l_n(\tilde{\psi})(\hat{\psi}_n-\psi^*)
$$
where $\tilde{\psi}=t\hat{\psi}_n+(1-t)\psi^*$. Normalizing  this expansion with $\sqrt{n}$ we obtain
$$
\sqrt{n}(\hat{\psi}_n-\psi^*)= -(\nabla_\psi^2l_n(\tilde{\psi}))^{-1} (\nabla_\psi l_n(\psi^*))\sqrt{n}.
$$

In order to obtain the asymptotic normality of the
maximum likelihood estimator we assume that exist an open 
neighborhood $B_r(\psi^*)$ of $\psi^*$ such that the following statements  hold.

\begin{itemize}
\item[H1] The functions $\psi\to A$ and $\psi\to p_\psi(Y_1|Y_0,i)$ are twice continuously differentiable
 on $B_r(\psi^*)$.
\item[H2] There exist functions $f_0,f_1,f_2$ such that 
 $$
 \sup_{\psi\in B_r(\psi^*)}\|\nabla_\psi p_\psi(y_1|y_0,i)\|\leq f_0(y_1,y_0),\ \
\sup_{\psi\in B_r(\psi^*)}\|\nabla_\psi^2 p_\psi(y_1|y_0,i)\|\leq f_1(y_1,y_0), 
$$ 
 
and 
$$
\sup_{\psi\in B_r(\psi^*)}\|\nabla_\psi p_\psi(y_1|y_0,i)\|\leq f_2(y_1,y_0),
$$ 
with $\mathbb{E}(f_s(Y_1,Y_0))<\infty$, $s=0,1$ and $\mathbb{E}(f_2(Y_1,Y_0)^{r})<\infty$, $r>2$.
\end{itemize}

\begin{teor} Under assumptions of Theorem 2.1 and H1-H2,
assume that  $J(\psi^*)=\mbox{var}(\nabla_\psi l(\psi^*))$
is non-singular and  
$\psi^*\in\mathring{\Psi}$. Then, as $n\to\infty$,
\begin{itemize}
\item[i)] $-(\nabla_\psi^2n^{-1}l_n(\tilde{\psi}))\to J(\psi^*)$, in probability. 
\item[ii)] $\sqrt{n}\nabla_\psi n^{-1} l_n(\psi^*)\to N(0,J(\psi^*))$, in  distribution.
\end{itemize}
Moreover, we conclude that $\sqrt{n}(\hat{\psi}_n-\psi^*)\to{\mathcal N}(0,J(\psi^*)^{-1})$, in distribution.
\end{teor}

\noindent{\bf Proof:}  Under H1-H2 if we take
$\varphi()=\frac{\partial^2 \log n^{-1}l_n}{\partial \psi^2}$ the Lemma \ref{LeyesMixing} implies 
that  
$$\sqrt{n}\nabla_\psi n^{-1}l_n(\psi^*)\to N(0,J(\psi^*))$$
and if $\varphi()=\frac{\partial n^{-1}\log l_n}{\partial \psi}$,
then 
$$-(\nabla_\psi^2n^{-1}l_n(\psi^*))\to J(\psi^*),\  a.s.$$

For a sequence $\psi_n\to\psi^*$  we can prove  
$$
\lim_{n\to\infty} \frac{1}{n}\nabla_\psi^2l_n({\psi}_n)-\frac{1}{n}\nabla_\psi^2l_n({\psi}_*)=0,
$$
in probability.\\

Let us first observe that 
$\frac{1}{n}\nabla_\psi^2l_n({\psi}_n)=\frac{1}{n}\sum_{k=1}^n\nabla_\psi^2 \log p_\psi(Y_k|Y_{0:k-1})$.
Another hand,
$$
\frac{\partial^2\log p_\psi}{\partial \psi_j\partial\psi_i}=\frac{1}{p_\psi}\frac{\partial^2 p_\psi}{\partial \psi_j\partial\psi_i}
-\frac{\partial\log p_\psi}{\partial \psi_j}\frac{\partial\log p_\psi}{\partial \psi_j}\frac{1}{p_\psi^2}
$$
hence,
\begin{eqnarray*}
\lefteqn{\frac{\partial^2\log p_{\psi_n}}{\partial \psi_j\partial\psi_i}-\frac{\partial^2\log p_{\psi^*}}{\partial \psi_j\partial\psi_i}}\\
&=&\left(\frac{1}{p_{\psi_n}}\frac{\partial^2 p_{\psi_n}}{\partial \psi_j\partial\psi_i}-\frac{1}{p_{\psi^*}}\frac{\partial^2 p_{\psi^*}}{\partial \psi_j\partial\psi_i}\right)+\left(\frac{\partial p_{\psi^*}}{\partial \psi_j}\frac{\partial p_{\psi^*}}{\partial \psi_j}\frac{1}{p_{\psi^*}^2}
-\frac{\partial p_{\psi_n}}{\partial \psi_j}\frac{\partial p_{\psi_n}}{\partial \psi_j}\frac{1}{p_{\psi_n}^2}\right)\\
&=& T_1+T_2.
\end{eqnarray*}

For term $T_2$, by definition of $\psi^*$,  
$\frac{\partial p_{\psi^*}}{\partial \psi_j}=0$ and by the Ergodic theorem we have
$$
\lim_{n\to\infty}\frac{1}{n}\sum_{k=1}^n\frac{\partial p_{\psi_n}}{\partial \psi_j}=0.
$$

For the term $T_1$,
$$
\frac{1}{p_{\psi_n}}\frac{\partial^2 p_{\psi_n}}{\partial \psi_j\partial\psi_i}-\frac{1}{p_{\psi^*}}\frac{\partial^2 p_{\psi^*}}{\partial \psi_j\partial\psi_i}=\frac{1}{p_{\psi_n}}\left(\frac{\partial^2 p_{\psi_n}}{\partial \psi_j\partial\psi_i}-\frac{\partial^2 p_{\psi^*}}{\partial \psi_j\partial\psi_i}\right)+\left(\frac{1}{p_{\psi_n}}-\frac{1}{p_{\psi^*}}\right)\frac{\partial^2 p_{\psi^*}}{\partial \psi_j\partial\psi_i} 
$$
Under equicontinuity of the sequence $\{p_{\psi_n}\}_{n\geq 1}$ we have, $p_{\psi_n}\to p_{\psi_{*}}$ and by conditions E1 and E7 
$\frac{1}{p_{\psi_n}}\to\frac{1}{p_{\psi^*}}$, a.s. Using H2 we obtain
$$
\mathbb{E}\left( \frac{\partial^2 p_{\psi^*}}{\partial \psi_j\partial\psi_i}\right)<\infty.
$$

Let 
$$
w(r,Y_{0:k})=\sup_{\psi_n\in B_r(\psi^*)}\left\lvert\frac{\partial^2 p_{\psi_n}(Y_k|Y_{0:k})}{\partial \psi_j\partial\psi_i}-
\frac{\partial^2 p_{\psi^*}(Y_k|Y_{0:k})}{\partial \psi_j\partial\psi_i}\right\rvert,
$$ 
proceeding as in Lemma 3 of Vandekerkhove \cite{vande}, by Markov inequality
\begin{eqnarray}
\label{InequalityMarkov}
\nonumber
\lefteqn{\mathbb{P}\left(\left\lvert\frac{1}{k}\sum_{k=1}^n\frac{\partial^2 p_{\psi_n}(Y_k|Y_{0:k})}{\partial \psi_j\partial\psi_i}-
\frac{1}{k}\sum_{k=1}^n\frac{\partial^2 p_{\psi^*}(Y_k|Y_{0:k})}{\partial \psi_j\partial\psi_i}\right\rvert
>\epsilon
\right)}\\
\nonumber
&\leq&\mathbb{P}\left(\frac{1}{k}\sum_{k=1}^nw(r,Y_{0:k})>\epsilon-\mathbb{E}(w(r,Y_{0:k}))
\right)+\mathbb{P}\left(\psi_n\not\in B_r(\psi^*)\right)\\
&\leq&\frac{\mathbb{E}(w(r,Y_{0:k}))}{\epsilon-\mathbb{E}(w(r,Y_{0:k}))}+\mathbb{P}\left(\psi_n\not\in B_r(\psi^*)\right).
\end{eqnarray}

The condicion H2  implies that $\mathbb{E}(w(r,Y_{0:k}))\leq 2f_1$. Using the Lebesgue continuity theorem,
we obtain that $\mathbb{E}(w(r,Y_{0:k}))\to 0$, as $n\to\infty$. The second term goes to $0$ as $n$ to infinity by strong convergence
of $\{\psi_n\}$ to $\psi^*$.
Hence (\ref{InequalityMarkov}) goes to $0$.\\
 
Finally, as $
\sqrt{n}(\hat{\psi}_n-\psi^*)=-\left(\nabla_\psi^2l_n(\tilde{\psi}))^{-1}(\nabla_\psi l_n(\psi^*)\right)\sqrt{n},
$
using i) the first factor in the above expression tends to  $J(\psi^*)$. The second factor converges 
weakly to $N(0,J(\psi^*))$ by ii). Slutsky's theorem  implies that 
$\sqrt{n}(\hat{\psi}_n-\psi^*)\to{\mathcal N}(0,J(\psi^*)^{-1})$.\cqd

\begin{ejemplo}(MS-AR gaussian again) 
 \end{ejemplo}
We employ the asymptotic results obtained to verify the validity of a likelihood test
for  identifing when the parameter $\rho$, of a MS-AR  is the zero vector. 
In this case the MS-AR process is a hidden Markov model.\\

Expanding $l_n(\rho)$ in Taylor series around $\hat{\rho}$, we have
$$
-2(l_n(\hat{\rho})-l_n(0))=\hat{\rho}^2\left(-\frac{\partial^2 l_n(\tilde{\rho})}{\partial \rho^2}\right)
$$
and by Theorem 2.1 $\hat{\rho}\sqrt{J(0)}\to\mathcal{N}(0,1)$ and
as $J(\tilde{\rho})/J(0)\to 1$ then
$\hat{\rho}^2{J(0)}\to \chi^2_1$.\\

\noindent{\bf Acknowledgments.} 
The author is grateful for the partial support on the projects Anillo ACT1112 and GEMINI-CONICYT 2012
NO. 32120025 CR 211291055 REXE 04464-14. Research facilities and hospitality in CIMFAV of the Universidad de Valpara\'{\i}so 
and by the sa\-bba\-tical support from the Universidad de Carabobo. The author also
thanks Dr K. Bertin for carefully reading a preliminary version.

\bibliographystyle{plain}

\begin{thebibliography}{10}

\bibitem{Andrews}
D.~Andrews.
\newblock {Non-strong mixing autoregressive procesess}.
\newblock {\em J. Appl Prob.}, 21:930--934, 1984.

\bibitem{LibCappe}
O.~Cappe, E.~Moulines, and T.~Ryd\'en.
\newblock {\em {\it Inference in Hidden Markov Models}}.
\newblock Springer-Verlag, 2005.

\bibitem{douc}
R.~Douc, E.~Moulines, and T.~Ryd{\'e}n.
\newblock {Asymptotic properties of the maximum likelihood estimator in
  autoregressive models with Markov regime}.
\newblock {\em Ann. Statist.}, 32:2254--2304, 2004.

\bibitem{dmixing}
P.~Doukhan.
\newblock {\em {\it Mixing: Propierties and Examples.}}, volume~85.
\newblock Lecture Notes in Statist., 1994.

\bibitem{Goldfeld}
S.~M. Goldfeld and R.~Quandt.
\newblock {A Markov Model for Switching Regressions}.
\newblock {\em Journal of Econometrics}, 1:3--16, 1973.

\bibitem{Hamilton}
J.D. Hamilton.
\newblock {A new approach to the economic analysis of non stationary time
  series and the business cycle}.
\newblock {\em Econometrica}, pages 357--384, 1989.

\bibitem{Kris}
V.~Krishnamurthy.
\newblock {Recursive Algorithms for estimation of hidden Markov Models with
  markov regime.}
\newblock {\em IEEE Trans. Information theory}, 48(2):458--476, 2002.

\bibitem{kris-Ryden}
V.~Krishnamurthy and T.~Ryd\'en.
\newblock {Consistent estimation of linear and non-linear autoregressive models
  with Markov regime.}
\newblock {\em Journal of Time Series Analysis}, 19:291--307, 1998.

\bibitem{luis3}
Ferm\'{\i}n L., R\'{\i}os, and L.~A. Rodr\'{\i}guez.
\newblock A Robbins Monro algorithm for nonparametric estimation of NAR process with Markov-Switching: consistency.
\newblock arXiv:1407.3747v6, 2014.

\bibitem{Rio1}
E.~Rio.
\newblock {\em {\it Th\'eorie asymptotique des processus faiblement
  d\'ependents}}, volume~31.
\newblock Springer-SMAI: Paris., 2000.

\bibitem{luis1}
R.~R\'{\i}os and L.~A. Rodr\'{\i}guez.
\newblock Penalized estimate of the number of states in gaussian linear ar with
  markov regime.
\newblock {\em Electronic Journal of Statistics}, pages 1111--1128, 2008.

\bibitem{ryn}
J.~Rynkiewicz.
\newblock {\em {\it Mod\'eles hybrides int\'egrant des r\'eseaux de neurones
  artificiels \`a des modeles de cha\^{i}nes de Markov cachee: application \`a
  la prediction de series temporelles }}.
\newblock PhD thesis, Universite Par{\'\i}s I, 2000.

\bibitem{Ramon}
R.~v.~Handel.
\newblock {\it Hidden Markov Models}.
\newblock Lecture notes: {\it https://www.princeton.edu/~rvan/}, 2008.

\bibitem{vande}
P.~Vandekerkhove.
\newblock {Consistent and asymptotically normal parameter estimates for hidden
  Markov mixtures of Markov models}.
\newblock {\em Bernoulli}, 11:103--129, 2005.

\bibitem{Yao}
J.~Yao and J.~G. Attali.
\newblock {On stability of nonlinear AR process with Markov switching}.
\newblock {\em Adv. Applied Probab}, 1999.

\end{thebibliography}

\end{document}